\documentclass[12pt]{article}
\usepackage{amssymb}
\newtheorem{theorem}{Theorem}

\newtheorem{proposition}[theorem]{Proposition}
\newtheorem{lemma}[theorem]{Lemma}
\newtheorem{corollary}[theorem]{Corollary}
\newtheorem{remark}[theorem]{Remark}
\newtheorem{remarks}[theorem]{Remarks}

\newcommand{\R}{\mathbb{R}}

\newcommand{\Sf}{\mathbb{S}}

\newcommand{\spa}{\mbox{span}}
\newcommand{\hess}{\mbox{Hess\,}}

\newcommand{\po}{{\hspace*{-1ex}}{\bf .  }}
\newcommand{\ii}{isometric immersion }
\newcommand{\iis}{isometric immersions }

\def\Z{{\mathbb Z}}
\def\<{\langle}

\def\>{\rangle}
\def\a{\alpha}

\def\bea{\begin{eqnarray*} }
\def\eea{\end{eqnarray*} }
\def\be{\begin{equation} }
\def\ee{\end{equation} }

\def\nap{\nabla^\perp}
\def\proof{\noindent{\it Proof: }}
\def\qed{\ifhmode\unskip\nobreak\fi\ifmmode\ifinner
\else\hskip5 pt \fi\fi\hbox{\hskip5 pt \vrule width4 pt
height6 pt  depth1.5 pt \hskip 1pt }}
\begin{document}
\title{\bigskip
\bigskip
Polar actions on compact Euclidean hypersurfaces.}
\author{Ion Moutinho $\&$ Ruy Tojeiro}
\date{}
\maketitle \noindent {\bf Abstract:} {\small Given an isometric
immersion $f\colon\,M^n\to \R^{n+1}$ of a compact Riemannian
manifold of dimension $n\geq 3$ into Euclidean space of dimension
$n+1$, we prove that the identity component $Iso^0(M^n)$ of the
isometry group $Iso(M^n)$ of $M^n$ admits an orthogonal
representation $\Phi\colon\,Iso^0(M^n)\to SO(n+1)$ such that
$f\circ g=\Phi(g)\circ f$ for every $g\in Iso^0(M^n)$. If $G$ is a
closed connected subgroup of $Iso(M^n)$ acting locally polarly on
$M^n$, we prove that $\Phi(G)$ acts polarly on $\R^{n+1}$, and we
obtain that $f(M^n)$ is given as $\Phi(G)(L)$, where $L$ is a
hypersurface of a section  which is invariant under the Weyl group
of the $\Phi(G)$-action. We also find several sufficient
conditions for such an $f$ to be a rotation hypersurface. Finally,
we show that compact Euclidean rotation hypersurfaces of
dimension $n\geq 3$ are characterized by their underlying warped product structure. } \vspace{3ex}\vspace{3ex}


\section[Introduction]{
Introduction}

   Let $G$ be a  connected
subgroup of the isometry group $Iso(M^n)$ of a compact Riemannian
manifold $M^n$ of dimension $n\geq 3$, which we always assume to be connected.   Given an isometric
immersion $f\colon\,M^n\to \R^{N}$ into Euclidean space of
dimension $N$, in general one can not expect $G$ to be realizable
as a group of rigid  motions of $\R^{N}$ that leave $f(M^n)$
invariant. Nevertheless, a fundamental fact for us is that in
codimension $N-n=1$ this is indeed the case.

\begin{theorem} \label{thm:main} Let $f\colon\,M^n\to \R^{n+1}$, $n\geq
3$, be a compact hypersurface. Then the identity component
$Iso^0(M^n)$ of the isometry group  of $M^n$ admits an orthogonal
representation $\Phi\colon\,Iso^0(M^n)\to SO(n+1)$ such that
$f\circ g=\Phi(g)\circ f$ for all $g\in Iso^0(M^n)$.
\end{theorem}

  Theorem \ref{thm:main}  may be regarded as a generalization of a classical result of
  Kobayashi \cite{ko}, who proved that a compact homogeneous
   hypersurface of Euclidean space must be a round sphere. In fact,
   the crucial step in the proof of Kobayashi's theorem is to show
   that the isometry group of the hypersurface can be realized as
   a closed subgroup of \mbox{$O(n+1)$}. The idea of the proof of Theorem \ref{thm:main}
   actually appears already in \cite{mpst}, where Euclidean $G$-hypersurfaces of cohomogeneity
   one, i.e., with principal
   orbits of codimension one, are considered.

    We  apply Theorem \ref{thm:main} to study compact Euclidean hypersurfaces
    $f\colon\, M^n\to \R^{n+1}$, $n\geq 3$, on which a connected closed subgroup
   $G$ of $Iso(M^n)$ acts locally polarly.
    Recall that an isometric
action of a compact Lie group $G$ on a Riemannian manifold $M^n$
is said to be  \emph{locally polar} if the distribution of normal
spaces to  principal orbits on the regular part of $M^n$ is
integrable. If $M^n$ is complete, this implies the existence of a
  connected complete immersed
submanifold $\Sigma$ of $M^n$ that intersects orthogonally all
$G$-orbits (cf. \cite{hlo}). Such a submanifold is called a
\emph{section}, and it is always a totally geodesic submanifold of
$M^n$. In particular, any isometric action of cohomogeneity one is
locally polar. The action is said to be {\em polar\/} if there
exists a closed and embedded section. Clearly, for orthogonal
representations there is no distinction between polar and locally
polar actions, for in this case sections are just affine
subspaces.

    It was shown in  \cite{bco}, Proposition $3.2.9$ that if a
    closed subgroup of $SO(N)$ acts polarly on $\R^{N}$ and leaves invariant a
    submanifold $f\colon\,M^n\to \R^{N}$,
    then its restricted action on $M^n$ is locally polar. Our
    next result states that any locally polar isometric action of a compact connected Lie group
    on a compact Euclidean
    hypersurface of dimension $n\geq 3$ arises in this way.

\begin{theorem} \label{cor:polar} Let $f\colon\,M^n\to \R^{n+1}$, $n\geq 3$,
be a compact hypersurface and let $G$ be a closed connected
subgroup of $Iso(M^n)$ acting locally polarly on $M^n$ with
cohomogeneity~$k$. Then there exists an orthogonal representation
$\Psi\colon\,G\to SO(n+1)$ such that $\Psi(G)$ acts polarly on
$\R^{n+1}$ with cohomogeneity $k+1$ and $f\circ g=\Psi(g)\circ f$
for every $g\in G$.    \end{theorem}

A natural problem that emerges is how to explicitly construct all
compact hypersurfaces $f\colon\,M^n\to \R^{n+1}$ that are invariant under
 a polar action of a closed subgroup  $G\subset SO(n+1)$. This is accomplished by the
following result. Recall that the Weyl group of the $G$-action  is defined as
$W=N(\Sigma)/Z(\Sigma)$, where  $\Sigma$ is a section, $N(\Sigma)=
\{g\in G\,|\,g\Sigma=\Sigma\}$ and $Z(\Sigma)=\{g\in G\,|\, gp=p,\,\forall\,p\in\Sigma\}$
is the intersection of the isotropy subgroups $G_p$, $p\in \Sigma$.

\begin{theorem} \label{thm:weyl} Let $G\subset SO(n+1)$ be a closed subgroup that acts
polarly on $\R^{n+1}$, let $\Sigma$ be a section and let $L$ be a
compact  immersed hypersurface  of $\Sigma$ which is invariant
under the Weyl group of the $G$-action. Then $G(L)$ is a compact
$G$-invariant immersed hypersurface of $\R^{n+1}$. Conversely, any
compact hypersurface $f\colon\,M^n\to \R^{n+1}$ that is  invariant
under a polar action of a closed subgroup of $SO(n+1)$ can be
constructed in this way.
\end{theorem}

   The simplest examples of hypersurfaces $f\colon\,M^n\to \R^{n+1}$ that are invariant under
   a polar action on $\R^{n+1}$ of a closed  subgroup of $SO(n+1)$ are the rotation
   hypersurfaces. These are invariant under a polar representation which is the sum of a
   trivial representation on a subspace $\R^k\subset \R^{n+1}$ and one which acts
   transitively on the unit sphere of its orthogonal complement $\R^{n-k+1}$.
   In this case any subspace $\R^{k+1}$ containing $\R^k$ is a section, and the
   Weyl group of the action is just the group $\Z^2$ generated by the reflection
   with respect to the ``axis"  $\R^k$. Thus, the condition that a hypersurface
   $L^k\subset \R^{k+1}$ be invariant under the Weyl group reduces in this case
   to  $L^k$ being symmetric with respect to $\R^k$.  We now give several sufficient
   conditions for a compact Euclidean hypersurface as in Theorem \ref{cor:polar}
   to be a  rotation hypersurface.

\begin{corollary} \label{cor:rev} Under the assumptions of Theorem
\ref{cor:polar}, any of the following additional conditions
implies that  $f$ is a rotation hypersurface:
\begin{itemize}
\item[$(i)$] there exists a totally geodesic (in $M^n$)
$G$-principal orbit; \item[$(ii)$] $k=n-1$; \item[$(iii)$] the
$G$-principal orbits are umbilical in $M^n$; \item[$(iv)$] there
exists a $G$-principal orbit with nonzero constant sectional
curvatures; \item[$(v)$] there exists a $G$-principal orbit with
positive sectional curvatures.
\end{itemize}
Moreover, $G$ is isomorphic to one of the closed subgroups of
$SO(n-k+1)$ that act transitively on $\Sf^{n-k}$.
\end{corollary}

    For a list of all closed subgroups of $SO(n)$ that act transitively
    on the sphere, see, e.g.,  \cite{eh}, p. 392.
    Corollary \ref{cor:rev} generalizes similar results in  \cite{ps},
     \cite{amn} and \cite{mpst} for compact Euclidean $G$-hypersurfaces of  cohomogeneity one.
      In part $(v)$, weakening the assumption to {\em non-negativity\/} of
      the sectional curvatures of some principal $G$-orbit implies $f$ to be a
      {\em multi-rotational hypersurface\/} in the sense of \cite{dn}:

     \begin{corollary} \label{cor:nneg} Under the assumptions of
     Theorem \ref{cor:polar}, suppose further that there exists a $G$-principal
orbit $Gp$ with nonnegative sectional curvatures. Then there exist
an orthogonal decomposition $\R^{n+1}=\bigoplus_{i=0}^k \R^{n_i}$
into $\tilde{G}$-invariant subspaces, where $\tilde{G}=\Psi(G)$,
and connected Lie subgroups $G_1,\ldots, G_k$ of $\tilde{G}$ such
that $G_i$ acts on $\R^{n_i}$, the action being transitive on
$\Sf^{n_i-1}\subset \R^{n_i}$, and the action of
$\bar{G}=G_1\times\ldots\times G_k$
 on $\R^{n+1}$ given by
 $$ (g_1\ldots
 g_k)(v_0,v_1,\ldots,v_k)=(v_0,g_1v_1,\ldots,g_kv_k)$$
 is orbit equivalent to the action of $\tilde{G}$.
In particular, if
 $Gp$ is flat then $n_i=2$ and $G_i$ is isomorphic to $SO(2)$ for $i=1,\ldots, k$.
\end{corollary}

    Finally, we apply some of the previous results to a problem  that at first sight has no
relation to isometric actions whatsoever.

     Let $f\colon\,M^n\to \R^{n+1}$ be a rotation hypersurface as described in the paragraph
     following Theorem \ref{thm:weyl}.
     Then the open and dense subset of $M^n$ that is
     mapped by $f$ onto the complement of the axis $\R^k$ is isometric to
     the warped product $L^k\times_{\rho} N^{n-k}$, where
     $N^{n-k}$ is the orbit of some fixed point
     $f(p)\in L^k$ under the action of $\tilde{G}$, and   the  warping function
     $\rho\,\colon\,L^k\to \R_+$  is a constant multiple of the  distance
     to  $\R^k$. Recall that a {\em warped product\/} $N_1\times_{\rho}N_2$ of Riemannian
manifolds
     $(N_1,\<\,\;,\,\>_{N_1})$ and $(N_2,\<\,\;,\,\>_{N_2})$ with {\em warping function\/} $\rho\colon\,N_1\to
     \R_+$ is the product manifold $N_1\times N_2$ endowed with
     the metric $$\<\,\;,\,\>=\pi_1^*\<\,\;,\,\>_{N_1}+(\rho\circ \pi_1)^2\pi_2^*\<\,\;,\,\>_{N_2},$$ where
     $\pi_i\colon\,N_1\times N_2\to N_i$, $1\leq i\leq 2$, denote
     the canonical projections.

         We prove that, conversely,  compact Euclidean rotation hypersurfaces
         of dimension $n\geq 3$ are characterized
         by their warped product structure.

    \begin{theorem} \label{thm:warp} Let $f\colon\, M^n\to \R^{n+1}$, $n\geq 3$,  be a
    compact hypersurface. If  there exists
     an isometry onto an open and dense subset $U\subset M^n$ of a warped product
     $L^k\times_{\rho}N^{n-k}$ with $N^{n-k}$ connected and complete
      (in particular if $M^n$ is isometric to a  warped product
$L^k\times_{\rho}N^{n-k}$)  then $f$ is a
rotation hypersurface.
\end{theorem}

     Theorem \ref{thm:warp} can be seen as a global  version in the
     hypersurface case of the local classification in \cite{dt} of isometric
     immersions in codimension $\ell\leq 2$ of warped products
     $L^k\times_{\rho}N^{n-k}$, $n-k\geq 2$, into Euclidean space.\vspace{2ex}

\noindent {\it Acknowledgment.\/} We are grateful to C. Gorodski
for helpful discussions.

\section[Proof of Theorem \ref{thm:main}]{Proof of Theorem
\ref{thm:main} }

    The proof of Theorem \ref{thm:main} relies on a result of Sacksteder
     \cite{sa} (see  also the proof in \cite{daj}, Theorem $6.14$, which is based on an
     unpublished manuscript by D. Ferus)
  according to which  a compact hypersurface  $f\colon\,M^n\to \R^{n+1}$, $n\geq 3$,
 is rigid whenever  the subset of totally geodesic points of $f$ does
not disconnect $M^n$.  Recall that $f$ is {\em rigid\/} if any
other isometric immersion $\tilde{f}\colon\, M^n\to \R^{n+1}$
differs from it by a rigid motion of $\R^{n+1}$. The proof of
Sacksteder's theorem actually shows more than the preceding
statement.
  Namely, let $f,\tilde{f}\colon\,M^n\to
\R^{n+1}$  be
 isometric immersions into $\R^{n+1}$ of a compact Riemannian manifold $M^n$,
$n\geq 3$, and let $\phi\colon\,T^\perp M^n_f\to T^\perp
M^n_{\tilde{f}}$ be the vector bundle isometry between the normal
bundles of $f$ and $\tilde{f}$ defined as follows. Given a unit
vector $\xi_x\in T^\perp_xM^n_f$ at $x\in M$, let $\phi(\xi_x)$ be
the unique unit vector in
 $T^\perp_xM^n_{\tilde{f}}$ such that
 $\{f_*X_1,\ldots, f_*X_n, \xi_x\}$ and $\{\tilde{f}_*X_1,\ldots, \tilde{f}_*X_n, \phi(\xi_x)\}$
  determine the same orientation in $\R^{n+1}$, where  $\{X_1,\ldots, X_n\}$ is any ordered basis
  of $T_xM^n$. Then it is shown  that
 \be\label{eq:sacks}\a_{\tilde{f}}(x)=\pm\phi(x)\circ \a_{f}(x)
\ee at each point $x\in M^n$, where $\a_{f}$ and $\a_{\tilde{f}}$
denote the second fundamental forms of $f$ and $\tilde{f}$,
respectively, with values in the normal bundle. The proof is based
on a careful analysis of the distributions on $M^n$ determined by
the {\em relative nullity subspaces\/}
$\Delta(x)=\mbox{ker}\;\a_{f}(x)$ and
$\tilde{\Delta}(x)=\mbox{ker}\;\a_{\tilde{f}}(x)$ of  $f$ and
$\tilde{f}$, respectively.  Rigidity of $f$ under the assumption
that the subset of totally geodesic points of $f$ does not
disconnect $M^n$ then follows immediately from the Fundamental
Theorem of Hypersurfaces (cf. \cite{daj}, Theorem
$1.1$).\vspace{1ex}

\noindent {\em Proof of Theorem \ref{thm:main}.} Given $g\in
Iso^0(M^n)$, let $\alpha_{f\circ g}$ denote the second fundamental
form of $f\circ g$. We claim that \be\label{eq:phig}\alpha_{f\circ
g}(x)=\phi_g(x)\circ \alpha_f(x)\ee for every $g\in Iso^0(M^n)$
and $x\in M^n$, where $\phi_g$ denotes the vector bundle isometry
between $T^\perp M^n_f$ and $T^\perp M^n_{{f}\circ g}$ defined as
in the preceding paragraph. On one hand,
\be\label{eq:sffs}\alpha_{f\circ
g}(x)(X,Y)=\alpha_f(gx)(g_*X,g_*Y)\ee for every $g\in Iso^0(M^n)$,
$x\in M^n$ and $X, Y\in T_xM^n$. In particular, this implies that
for any  fixed $x\in M^n$ the map $\Theta_x\colon\,Iso^0(M^n)\to
\mbox{Sym}(T_xM^n\times T_xM^n\to T^\perp_xM^n_f)$ into the vector
space of symmetric bilinear maps of $T_xM^n\times T_xM^n$ into
$T^\perp_xM^n_f$, given by
$$\Theta_x(g)(X,Y)=\phi_g(x)^{-1}(\alpha_{f\circ
g}(x)(X,Y))=\phi_g(x)^{-1}(\alpha_f(gx)(g_*X,g_*Y))$$ for any
$X,Y\in T_xM^n$, is continuous. On the other hand, by the
preceding remarks on Sacksteder's theorem,  either $\a_{f\circ
g}(x)= \phi_g(x)\circ \a_{f}(x)$  or $\a_{f\circ g}(x)=
-\phi_g(x)\circ \a_{f}(x)$. Thus $\Theta_x$ is a continuous map
taking values in $\{\a_f(x), -\a_f(x)\}$, hence it must be
constant because $Iso^0(M^n)$ is connected. Since
$\Theta_x(\mbox{id})=\a_f(x)$, our claim follows.

    We conclude that for each $g\in Iso^0(M^n)$ there exists a rigid motion
    $\tilde{g}\in  Iso(\R^{n+1})$
    such that $f\circ g= \tilde{g}\circ f$. It now follows from
    standard arguments that $g\mapsto \tilde{g}$ defines a
    Lie-group homomorphism $\Phi\colon\,Iso^0(M^n)\to
    Iso(\R^{n+1})$, whose image must lie in $SO(n+1)$ because it
    is compact and connected.\qed

     \begin{remarks}{\em $(i)$ Theorem \ref{thm:main} is also true for compact hypersurfaces
of dimension  $n\geq 3$ of hyperbolic space, as well as for {\em
complete\/} hypersurfaces of dimension $n\geq 4$ of the sphere. It
also holds for  complete   hypersurfaces of dimension
$n\geq 3$ of both Euclidean and hyperbolic spaces, under the
additional assumption that they do not carry a complete leaf of
dimension $(n-1)$ or $(n-2)$ of their relative nullity
distributions. In fact, the proof of Theorem \ref{thm:main}
carries over in exactly the same way for these cases, because so
does equation (\ref{eq:sacks}) (cf. \cite{daj}, p.
96-100).\vspace{1ex}\\
 $(ii)$   Clearly, Theorem \ref{thm:main} does not hold
    for \iis $f\colon\, M^n\to \R^{n+\ell}$ of codimension
    $\ell\geq 2$. Namely, counterexamples can be easily constructed, for instance, by  means of    compositions $f=h\circ g$ of isometric immersions $g\colon\, M^n\to \R^{n+1}$
    and $h\colon\, V\to \R^{n+\ell}$, with
    $V\subset \R^{n+1}$ an open subset containing $g(M^n)$.}
    \end{remarks}

     The next result gives a sufficient condition for
     rigidity  of a compact hypersurface $f\colon\,M^n\to \R^{n+1}$ as in Theorem \ref{thm:main} in terms of
$\tilde{G}=\Phi(Iso^0(M^n))$.

\begin{proposition} \label{prop:rig} Under the assumptions of Theorem
\ref{thm:main}, suppose  that $\tilde{G}=\Phi(Iso^0(M^n))$ does
not have a fixed vector. Then $f$ is free of totally geodesic
points. In particular, $f$ is rigid.
\end{proposition}
\proof  Suppose that the subset $B$ of totally geodesic points of
$f$ is nonempty. Since $\alpha_{f\circ g}(x)=\phi_g(x)\circ
\alpha_f(x)$ for every $g\in G=Iso^0(M^n)$ and  $x\in M^n$ by
(\ref{eq:phig}), $B$ coincides with the set of totally geodesic
points of $f\circ g$ for every $g\in G$. In view of
(\ref{eq:sffs}), this amounts to saying that $B$ is $G$-invariant.
Thus, if $Gp$ is the orbit of a point $p\in B$ then $Gp\subset B$.
Since  $Gp$ is connected, it follows from  \cite{dg}, Lemma 3.14
that $f(Gp)$ is contained in a hyperplane ${\cal H}$ that is
tangent to $f$ along $Gp$. Therefore, a unit vector $v$ orthogonal
to ${\cal H}$ spans $T_{gp}^\perp M^n_f$ for every $gp\in Gp$.
Since $\tilde{g}T_{p}^\perp M^n_f=T_{gp}^\perp M^n_f$ for every
$\tilde{g}=\Phi(g)\in \tilde{G}$, because $f$ is equivariant with
respect to $\Phi$, the connectedness of $\tilde{G}$ implies that
it must fix $v$. \vspace{2ex} \qed

    Proposition \ref{prop:rig} implies, for instance,  that if a closed
   connected subgroup of $Iso(M^n)$ acts on
   $M^n$ with cohomogeneity one  then either $f$ is a rotation hypersurface over a plane curve
   or it is free
   of totally geodesic points, and in particular it is rigid (see \cite{mpst}, Theorem $1$).
   As another consequence we have:

\begin{corollary} \label{cor:rig} Let $f\colon\, M^3\to \R^{4}$ be a compact hypersurface.
If $f$ has a totally geodesic point (in particular, if it is not
rigid) then either $M^3$ has finite isometry group or $f$ is a
rotation hypersurface.
\end{corollary}
\proof Let $\Phi\colon\,Iso^0(M^3)\to SO(4)$ be the orthogonal
representation given by Theorem~\ref{thm:main}. By Proposition
\ref{prop:rig}, if $f$ has a totally geodesic point  then
$\tilde{G}=\Phi(Iso^0(M^3))$ has a fixed vector $v$, hence it can
be regarded as a subgroup of $SO(3)$. Therefore, either the
restricted action of $\tilde{G}$ on $\{v\}^\perp$ has also a fixed
vector or it is transitive on the sphere. In the first case,
either $Iso^0(M^3)$ is trivial, that is, $Iso(M^3)$ is finite, or
$\tilde{G}$ fixes a two dimensional subspace $\R^2$ of $\R^4$, in
which case $f$ is a rotation hypersurface over a surface in a
half-space $\R^3_+$ with $\R^2$ as boundary. In the latter case,
$f$ is a rotation hypersurface over a plane curve in a half-space
$\R^2_+$ having $\spa\{v\}$ as boundary.\qed

\section[Proof of Theorem \ref{cor:polar}]{Proof of Theorem \ref{cor:polar}}

 For the proof of Theorem \ref{cor:polar}, we recall from
Theorem \ref{thm:main} that  there exists an orthogonal
representation $\Phi\colon\,Iso^0(M^n)\to SO(n+1)$ such that
$f\circ g=\Phi(g)\circ f$ for every $g\in Iso^0(M^n)$. Since $G$
is connected we have $G\subset Iso^0(M^n)$, hence it suffices to
prove that $\tilde{G}=\Phi(G)$ acts polarly on $\R^{n+1}$ with
cohomogeneity $k+1$ and then set $\Psi=\Phi|_G$.

 We claim that there exists a principal orbit
$Gp$ such that
     the position vector $f$ is nowhere tangent to $f(M^n)$ along
     $Gp$, that is, $f(g(p))\not\in f_*(g(p))T_{g(p)}M^n$ for any
     $g\in G$. In order to prove our claim we need the following observation.

     \begin{lemma} \label{le:tang} Let $f\colon\, M^n\to \R^{n+1}$ be a
     hypersurface. Assume that the position vector is tangent to
     $f(M^n)$ on an open subset $U\subset M^n$. Then the index of
     relative nullity $\nu_f(x)=\mbox{dim }\Delta_f(x)$ of $f$ is positive at any point
     $x\in U$.
     \end{lemma}
     \proof  Let $Z$ be a vector field on $U$ such that
     $f_*(p)Z(p)=f(p)$ for any $p\in U$ and let $\eta$ be a local unit normal vector field to $f$.
      Differentiating $\<\eta,f\>=0$ yields $\<AX,Z\>=0$ for any
      tangent vector field $X$ on $U$, where $A$ denotes the shape operator of $f$ with respect to $\eta$.
       Thus $AZ=0$.\vspace{1ex}\qed

     Going back to the proof of the claim, since $M^n$ is a compact Riemannian
     manifold isometrically immersed in
     Euclidean space as a hypersurface,  there exists an open subset
     $V\subset M^n$ where the sectional curvatures of $M^n$ are
     strictly positive. In particular, the index of relative nullity of $f$ vanishes
     at every $x\in V$.
     If the position vector were tangent to $f(M^n)$
     at every regular point of $V$, it would be tangent to $f(M^n)$ everywhere on $V$, because
     the set of regular points is dense on $M^n$. This is in contradiction with
     Lemma~\ref{le:tang} and shows that there must exist a regular
     point $p\in V$ such that $f(p)\not\in f_*(p)T_pM^n$. Since $f$
     is equivariant, we must have in fact that $f(g(p))\not\in
     f_*(g(p))T_{g(p)}M^n$ for any $g\in G$ and the  claim is proved.

       Now let  $Gp$ be a principal orbit such that
     the position vector $f$ is nowhere tangent to $f(M^n)$ along
     $Gp$. Then the normal bundle of $f|_{Gp}$ splits (possibly non-orthogonally) as
     $\mbox{span}\{f\}\oplus f_*T^\perp Gp$.
Let $\xi$ be an
     equivariant normal vector field to $Gp$ in $M^n$ and let
     $\eta=f_*\xi$. Then, denoting $\Phi(g)$
     by $\tilde{g}$ and identifying  $\tilde{g}$ with its
     derivative at any point, because it is linear, we have
     \be\label{eq:equiv}\tilde{g}\eta(p)=(f\circ
     g)_*(p)\xi(p)=f_*(gp)g_*(p)\xi(p)=f_*(gp)\xi(gp)=\eta(gp)\ee
     for any $g\in G$. In particular,
     $\<\eta(gp),f(gp)\>= \<\tilde{g}\eta(p),\tilde{g}f(p)\> =\< \eta(p),f(p)\> $
     for every $g\in G$, that is, $\<\eta,f\>$ is constant on $Gp$. It follows that
  \be\label{eq:const} X\<\eta,f\>=0, \,\,\,\mbox{for
     any}\,\,\,X\in TGp,\ee
     and hence
     \be\label{eq:comp1}\<\tilde{\nabla}_X
     \eta,f\>=X\<\eta,f\>-\<\xi,X\>=0,\ee
      where $\tilde{\nabla}$ denotes the derivative in $\R^{n+1}$. On the other hand,
      since $G$ acts locally polarly on $M^n$, we have that
     $\xi$ is parallel in the normal connection of $Gp$ in
     $M$. Therefore
     \be\label{eq:comp2}(\tilde{\nabla}_X \eta)_{{f}_{*}T^{\perp} Gp}=f_*({\nabla_X\xi})_{T^\perp
     Gp}=0,\ee
   where $\nabla$ is the Levi-Civita connection of $M^n$; here, writing a vector subbundle as
   a subscript
of a vector field indicates taking its orthogonal projection onto
that subbundle. It follows from (\ref{eq:comp1}) and
(\ref{eq:comp2}) that $\eta$ is parallel in the normal connection
of $f|_{Gp}$. On the other hand, the position vector $f$ is clearly also an equivariant normal
  vector  field of $f|_{Gp}$ which is parallel in the normal connection.

     Thus, we have shown that there exist equivariant normal vector fields to $\tilde{G}(f(p))=f({Gp})$ that form
     a basis of the normal spaces at each point and which are parallel in the normal connection. The statement
     now follows from the next known result, a  proof of which is included for completeness.
\begin{lemma} \label{le:pol1} Let $\tilde{G}\subset SO(n+1)$ have an orbit $\tilde{G}(q)$ along which
there exist equivariant normal vector fields that form a basis of
the normal spaces at each point and which are parallel in the
normal connection. Then $\tilde{G}$ acts polarly.
\end{lemma}
\proof Since there exist equivariant normal vector fields to
$\tilde{G}(q)$ that form a basis of the normal spaces at each
point, the isotropy group acts trivially on each normal space,
hence $\tilde{G}(q)$ is a principal orbit. We now show that the
normal space $T^\perp_q \tilde{G}(q)$ is a section, for which it
suffices to show that any Killing vector field $X$ induced by
$\tilde{G}$ is everywhere orthogonal to $T^\perp_q \tilde{G}(q)$.
Given $\xi_q\in T^{\perp}_q \tilde{G}(q)$, let $\xi$ be an
equivariant normal vector field to $\tilde{G}(q)$ extending
$\xi_q$, which is also parallel in the normal connection. Then,
denoting by $\phi_t$ the flow of $X$ and setting $c(t)=\phi_t(q)$
we have
$$\left(\frac{d}{dt}|_{t=0}\phi_t(\xi_q)\right)_{T^{\perp}_q \tilde{G}(q)}
=\left(\frac{D}{dt}|_{t=0}\xi(c(t))\right)_{T^{\perp}_q
\tilde{G}(q)}=\nabla^\perp_{X(q)}\xi=0,$$ where ${\displaystyle \frac{D}{dt}}$
denotes the covariant derivative in Euclidean space along $c(t)$.\qed

     \begin{remark}{\em  A closed subgroup $G\subset SO(n+\ell)$, $\ell\geq 2$, that acts
     non-polarly  on $\R^{n+\ell}$ may leave invariant a compact submanifold $f\colon\, M^n\to
     \R^{n+\ell}$  and induce a locally polar action on $M^n$. For instance, consider a compact
     submanifold $f\colon\, M^n\to  \R^{n+2}$ that is invariant by the action
     of a closed subgroup $G\subset SO(n+2)$ that acts
     non-polarly on $\R^{n+2}$ with cohomogeneity three. Then the induced action of $G$
     on $M^n$ has cohomogeneity one, whence  is locally polar. Moreover, taking
     $f$ as a compact hypersurface of $\Sf^{n+1}$ shows also that Theorem
     \ref{cor:polar} is no longer true if $\R^{n+1}$ is replaced by $\Sf^{n+1}$.
     }    \end{remark}

     Theorem \ref{cor:polar} yields the following obstruction for
   the existence of an isometric immersion in codimension one into
   Euclidean space of a compact Riemannian manifold acted on locally polarly
   by a closed connected Lie subgroup of its isometry group.

   \begin{corollary} \label{cor:excep} Let $M^n$ be a compact Riemannian manifold of
   dimension $n\geq 3$ acted on locally polarly by a closed connected  subgroup $G$ of
its isometry group. If $G$ has an exceptional orbit then $M^n$ can
not be isometrically immersed in Euclidean space as a
hypersurface.
   \end{corollary}
\proof  Let $f\colon\,M^n\to\R^{n+1}$ be an isometric immersion of
a compact Riemannian manifold acted on locally polarly by a closed
connected
 subgroup $G$ of its isometry group. We will prove that $G$ can
not have any exceptional orbit. By Theorem~\ref{cor:polar} there
exists an orthogonal representation $\Psi\colon\,G\to SO(n+1)$
such that $\tilde{G}=\Psi(G)$ acts polarly on $\R^{n+1}$ with
cohomogeneity $k+1$ and $f\circ g=\Psi(g)\circ f$ for every $g\in
G$. Let $Gp$ be a nonsingular orbit. Then $Gp$ has maximal
dimension among all $G$-orbits, and hence $\tilde{G}f(p)=f(Gp)$
has maximal dimension among all $\tilde{G}$-orbits. Since polar
representations are known to admit no exceptional orbits (cf.
\cite{bco}, Corollary $5.4.3$), it follows  that  $\tilde{G}f(p)$
is a principal orbit. Then, for any $g$ in the isotropy subgroup
$G_p$ we have that $\tilde{g}=\Psi(g)\in \tilde{G}_{f(p)}$, thus
for any $\xi_p\in T_p^\perp Gp$ we obtain
$$f_*(p)\xi_p=\tilde{g}_*f_*(p)\xi_p=(\tilde{g}\circ
f)_*(p)\xi_p=(f\circ
g)_*(p)\xi_p=f_*(gp)g_*\xi_p=f_*(p)g_*\xi_p.$$ Since $f_*(p)$ is
injective, then $g_*\xi_p=\xi_p$. This shows that the slice
representation, that is, the action of the isotropy group $G_p$ on the normal space
$T_p^\perp G(p)$  to the orbit $G(p)$ at $p$, is trivial. Thus $Gp$ is a
principal orbit. \vspace{2ex}\qed

\section[Isoparametric submanifolds]{Isoparametric
    submanifolds}

    We now recall some results on isoparametric
    submanifolds and derive a few additional facts on them that
    will be needed for the proofs of Theorem \ref{thm:weyl} and Corollaries~\ref{cor:rev}
    and \ref{cor:nneg}.

   Given an isometric immersion
$f\colon\, M^{n}\to \R^{N}$ with flat normal bundle, it is
well-known (cf. \cite{st}) that for each point $x \in M^{n}$ there
exist an integer $s=s(x) \in \{1, \ldots , n\}$ and a uniquely
determined subset $H_{x}=\{ \eta_{1}, \ldots , \eta_{s} \}$ of
$T_{x}^{\perp}M_f^n$ such that $T_{x}M^n$ is the orthogonal sum of
the nontrivial subspaces
  \[ E_{\eta_i}(x)=\{ X \in T_{x}M^n \, \colon\,\alpha(X,Y)=\langle X,Y \rangle
\, \eta_{i}, \,\, \mbox{ for all } \, Y \in T_{x}M^n \},\,\, 1\leq
i\leq s. \] Therefore, the second fundamental form of $f$ has the
simple representation
\begin{equation}\label{eq:2forma}
 \alpha(X,Y)=\sum_{i=1}^{s} \langle X_i,Y_i \rangle \, \eta_{i},
\end{equation}
or equivalently, \begin{equation}\label{eq:2formab} A_\xi
X=\sum_{i=1}^s\<\xi,\eta_i\>X_i, \end{equation}
 where $X\mapsto X_i$ denotes
orthogonal projection onto $E_{\eta_i}$. Each $\eta_{i} \in H_{x}$
is called a {\em principal normal\/} of $f$ at $x$. The Gauss
equation takes the form
\begin{equation} \label{eq:gauss}
R(X,Y)=\sum_{i,j=1}^{s}\langle \eta_i ,\eta_j \rangle X_i \wedge
Y_j ,
\end{equation}
where  $(X_i\wedge Y_j)Z=\<Z,Y_j\>X_i-\<Z,X_i\>Y_j$.

\begin{lemma}\label{le:n1} Let $f\colon\,M^n\to \R^N$
be an  \ii with flat normal bundle of a Riemannian manifold with
constant sectional curvature $c$. Let $\eta_1,\ldots,\eta_s$ be
the distinct principal normals of $f$ at $x\in M^n$. Then
\begin{itemize} \item[$(i)$] There
exists at most one $\eta_i$ such that $|\eta_i|^2=c$.
\item[$(ii)$] For all $i,j,k\in \{1,\ldots, s\}$ with $i\neq j\neq
k\neq i$ the vectors $\eta_i-\eta_k$ and $\eta_j-\eta_k$ are
linearly independent. \end{itemize}
\end{lemma}
\proof It follows from (\ref{eq:gauss}) that $\mbox{$\langle
\eta_i ,\eta_j \rangle=c$}$   for all $i,j\in \{1,\ldots, s\}$
with $i\neq j$.  If $|\eta_i|^2=|\eta_j|^2=c$, this gives
$|\eta_i-\eta_j|^2=0$.
\vspace{1ex}\\
$(iii)$  Assume that there exist $\lambda\neq 0$ and $i,j,k\in
\{1,\ldots, s\}$ with $i\neq j\neq k\neq i$ such that
$\eta_i-\eta_k=\lambda(\eta_j-\eta_k)$. Then
$$|\eta_i|^2-c=\<\eta_i-\eta_k,\eta_i\>=\lambda\<\eta_j-\eta_k,\eta_i\>=0,$$
and similarly $|\eta_j|^2=c$, in contradiction with
$(i)$.\vspace{2ex}\qed

    Let $f\colon\,M^n\to \R^N$
be an  \ii with flat normal bundle and let $H_{x}=\{ \eta_{1},
\ldots , \eta_{s} \}$ be the set of principal normals of $f$ at
$x\in M^n$. If the map $M^{n} \to \{1, \ldots ,n \}$ given by $x
\mapsto \# H_{x}$ has a constant value $s$ on an open subset
$U\subset M^n$, there exist smooth normal vector fields $\eta_{1},
\ldots , \eta_{s}$ on $U$ such that $H_{x}= \{ \eta_{1}(x), \ldots
, \eta_{s}(x) \}$ for any $x\in U$. Furthermore, each
$E_{\eta_i}=(E_{\eta_i}(x))_{x \in U}$ is a $C^{\infty}$-
subbundle of $TU$ for $1 \leq i\leq s$. The following result is
contained in \cite{dn}, Lemma $2.3$.

\begin{lemma}\po \label{cor:regu0}  Let
$f\colon\,M^n\to\R^{n+p}$ be an \ii with flat normal bundle and a
constant number $s$ of principal normals $\eta_1,\ldots ,\eta_s$
everywhere. Assume that for a fixed $i\in \{1,\ldots, s\}$ all
principal normals $\eta_j$, $j\neq i$, are parallel in the normal
connection along $E_{\eta_i}$ and that  the vectors
$\eta_i-\eta_{j}$ and $\eta_i-\eta_\ell$ are everywhere pairwise
linearly independent for any pair of indices \mbox{$1\leq j\neq
\ell\leq s$} with $j, \ell\neq i$. Then $E_{\eta_i}^\perp$ is
totally geodesic.
\end{lemma}

\proof The Codazzi equation yields
 \be \label{eq:codz2}
\<\nabla_{X_\ell}X_j,X_i\> (\eta_i-\eta_j)=
\<\nabla_{X_\ell}X_j,X_i\> (\eta_i-\eta_\ell), \,\,\,\,i\neq j\neq
k\neq i,\ee and \be\label{eq:codz}
\nap_{X_i}\eta_j=\<\nabla_{X_j}X_j,X_i\>(\eta_j-\eta_i),
\;\;\;i\neq j, \ee for all unit vectors $X_i\in E_{\eta_i}$,
$X_j\in E_{\eta_j}$ and $X_{\ell}\in
E_{\eta_{\ell}}$.\vspace{2ex}\qed

An \ii $f\colon\,M^n\to \R^N$ is  called {\em isoparametric\/} if
it has flat normal bundle and the principal curvatures of $f$ with
respect to every parallel normal vector field along any curve in
$M^n$ are constant (with constant multiplicities). The following
facts on isoparametric submanifolds are due to Str\"ubing
\cite{st}.

\begin{theorem}\label{thm:s} Let $f\colon\,M^n\to \R^N$
be an  isoparametric isometric immersion.  Then
\begin{itemize} \item[$(i)$] The number of principal normals is
constant on $M^n$.
 \item[$(ii)$] The first normal spaces, i.e, the subspaces of the normal spaces spanned by the image of
 the second fundamental form, determine  a parallel subbundle
 of the normal bundle.
\item[$(iii)$] The subbundles $E_{\eta_i}$, $1\leq i\leq s$, are
 totally geodesic and the principal normals $\eta_1,\ldots,
 \eta_s$ are parallel in the normal connection.
\item[$(iv)$] The subbundles $E_{\eta_i}$, $1\leq i\leq s$, are
parallel if and only if $f$ has parallel second fundamental form.
\item[$(v)$] If the principal normals $\eta_1,\ldots,\eta_s$
satisfy $\langle \eta_i ,\eta_j \rangle\geq 0$ everywhere then $f$
has parallel second fundamental form and $\langle \eta_i ,\eta_j
\rangle= 0$ everywhere.
\end{itemize}
\end{theorem}

The next result will be used in the proofs of Corollaries
\ref{cor:rev} and \ref{cor:nneg}.

\begin{proposition} \label{prop:isop1} Let $f\colon\, M^n\to \R^{N}$, $n\geq 2$, be
a compact isoparametric submanifold.
\begin{itemize} \item[$(i)$] If $M^n$ has constant sectional
curvature $c$, then either  $c>0$ and $f(M^n)$ is a round sphere
or  $c=0$ and $f(M^n)$ is an extrinsic  product of circles.
\item[$(ii)$] If $M^n$ has nonnegative sectional curvatures, then
$f(M^n)$ is an extrinsic product of round spheres or circles. In
particular, if $M^n$ has positive sectional curvatures then
$f(M^n)$ is a round sphere.
     \end{itemize}
\end{proposition}
\proof  By Theorem \ref{thm:s}, the number of distinct principal
normals $\eta_1,\ldots, \eta_s$ of $f$ is constant on $M^n$ and
all of them are parallel in the normal connection. Moreover, the
subbundles $E_{\eta_i}$, $1\leq i\leq s$, are
 totally geodesic. If $M^n$ has constant sectional
curvature~$c$, then it follows from Lemmas \ref{le:n1} and
\ref{cor:regu0} that also $E_{\eta_i}^\perp$ is totally geodesic
for $1\leq i\leq s$, and hence the sectional curvatures along
planes spanned by vectors in different subbundles vanish.
Therefore $c=0$, unless $s=1$ and  $f$ is umbilic, in which case
$c>0$ and $f(M^n)$ is a round sphere. Furthermore, if $c=0$ and
$E_{\eta_i}$ has rank at least $2$ for some $1\leq i\leq s$ then
the sectional curvature along a plane tangent to $E_{\eta_i}$ is
$|\eta_i|^2=0$, in contradiction with the compactness of $M^n$.
Hence $E_{\eta_i}$ has rank $1$ for $1\leq i\leq s$. We conclude
that the universal covering of $M^n$ is isometric to $\R^n$, and
that $f\circ \pi$ splits as a product of circles by Moore's Lemma
\cite{mo}, where $\pi\colon\,\R^n\to M^n$ is the covering map.

  Assume now that $M^n$ has nonnegative sectional curvatures. It
  follows from (\ref{eq:gauss}) that $\langle
  \eta_i,\eta_j\rangle\geq 0$ for $1\leq i\neq j\leq s$, whence
  $f$ has parallel second fundamental form  and all subbundles $E_{\eta_i}$, $1\leq i\leq s$,
  are parallel by parts
  $(iv)$ and $(v)$ of Theorem  \ref{thm:s}. We obtain from the de Rham decomposition theorem
  that the
universal covering of $M^n$ splits isometrically as
$M_1^{n_1}\times \cdots\times M_s^{n_s}$, where each factor
$M_i^{n_i}$ is either $\R$ if $n_i=1$ or a sphere  $\Sf_i^{n_i}$
of curvature $|\eta_i|^2$ if $n_i\geq 2$. Moreover, if
$\pi\colon\,M_1^{n_1}\times \cdots\times M_s^{n_s}\to M^n$ denotes
the covering map, then Moore's Lemma implies that  $f\circ \pi$
splits as $f\circ \pi=f_1\times\cdots\times f_s$, where
$f_i(M_i^{n_i})$ is a round sphere or circle for $1\leq i\leq s$.
\vspace{2ex}\qed

To every compact isoparametric submanifold $f\colon\, M^n\to
\R^{N}$ one can associate a finite group, its Weyl group, as
follows. Let $\eta_1, \ldots, \eta_g$ denote the principal
 normal vector fields of $f$.
For $p\in M^n$, let $H_j(p)$, $1\leq j\leq g$,  be the focal hyperplane of $T^\perp_pM^n$
given by the equation $\<\eta_j(p),\,\,\>=1$. Then one can show that the reflection on the
affine normal space $p+T_p^\perp M^n$ with respect to each affine focal hyperplane $p+H_i(p)$
leaves
$\bigcup_{j=1}^g(p+H_j(p))$ invariant, and thus the set of all such reflections generate a
finite group, the Weyl group of $f$ at $p$. Moreover,  the Weyl groups of $f$ at different
points are conjugate by the parallel transport with respect to the normal connection, hence
a well-defined Weyl group $W$ can be associated to  $f$.  We refer to \cite{pt2} for details.
In the proof of Theorem \ref{thm:weyl} we will need the following property of the Weyl group of
an isoparametric submanifold.

\begin{proposition} \label{prop:weyl} Let $f\colon\, M^n\to \R^{N}$ be a compact
isoparametric submanifold and let $W(p)$ be its Weyl group at $p\in M^n$. Assume that $W(p)$
leaves invariant an affine hyperplane ${\cal H}$ orthogonal to $\xi\in T^\perp_pM$. Then
$f(M^n)$ is contained in the affine hyperplane of $\R^N$ through $p$ orthogonal to $\xi$.
\end{proposition}
\proof It follows from the assumption that ${\cal H}$ is orthogonal to every focal hyperplane
$p+H_j(p)$, $1\leq j\leq g$, of $f$ at $p$. For  $q\in {\cal H}$, let
$Q=\sum_{g\in W(p)} gq\in {\cal H}$. Then $Q$ is a fixed point of $W(p)$,
hence it lies in the intersection
$\bigcap_{j=1}^g (p+H_j(p))$
of all affine focal hyperplanes of $f$ at $p$. We obtain that the line through $Q$
orthogonal to ${\cal H}$ lies in $\bigcap_{j=1}^g (p+H_j(p))$.
Therefore $\<\eta_j,Q+\lambda\xi\>=1$ for every $\lambda\in \R$, $1\leq j\leq g$,
which implies that $\<\eta_j,\xi\>=0$, for every $1\leq j\leq g$. Now extend $\xi$
to a parallel vector field along $M^n$ with respect to the normal connection. Since
the principal normal vector fields $\eta_1, \ldots, \eta_g$ of $f$ are parallel with
respect to the normal connection
 by Theorem \ref{thm:s}-$(iii)$, it follows that $\<\eta_j,\xi\>=0$ everywhere,
 and hence the shape operator $A_\xi$ of $f$ with respect to $\xi$ is identically
 zero by (\ref{eq:2formab}). Then $\xi$ is constant in $\R^N$ and the conclusion
 follows.\vspace{2ex}\qed

A rich source of isoparametric submanifolds is provided by the
following result of Palais and Terng (see \cite{pt}, Theorem
$6.5$).

\begin{proposition} \label{prop:pt1} If a closed subgroup $G\subset SO(N)$ acts polarly on
$\R^{N}$  then any of its
principal orbits is an isoparametric submanifold of $\R^N$.
\end{proposition}

To conclude this section, we point out that if  $G\subset SO(N)$ acts polarly on $\R^{N}$
and $\Sigma$ is a section, then the the Weyl group $W=N(\Sigma)/Z(\Sigma)$ of the
$G$-action coincides with the Weyl group $W(p)$  just defined of any principal orbit
$Gp$, $p\in \Sigma$, as an isoparametric submanifold of $\R^N$ (cf. \cite{pt2}).

    \section[Proof of Theorem \ref{thm:weyl}]{Proof of Theorem \ref{thm:weyl}}

We first prove the converse. Let $G\subset SO(n+1)$ act polarly on $\R^{n+1}$,
let $\Sigma$ be a section of the $G$-action  and let  $M^n\subset \R^{n+1}$ be a
$G$-invariant immersed hypersurface. It suffices to prove that $\Sigma$ is transversal to $M^n$,
 for then $L=\Sigma\cap M^n$ is a compact hypersurface of $\Sigma$ that is invariant under
 the Weyl group $W$ of the action and $M^n=G(L)$.

   Assume, on the contrary, that transversality does not hold. Then there exists
   $p\in \Sigma\cap M^n$ with $T_p\Sigma\subset T_pM^n$.  Fix $v\in T_p\Sigma$ in
   a principal orbit of the slice representation at $p$ and
   let $\gamma\colon\,(-\epsilon, \epsilon)\to M^n$ be a smooth curve with $\gamma(0)=p$
   and  $\gamma'(0)=v$. Since $M^n$ is $G$-invariant, it contains
   $\{g\gamma(t)\,|\,g\in G_p, t\in (-\epsilon, \epsilon)\}$. Therefore,
   $T_pM^n$ contains $G_pv$, and hence $\R v \oplus T_vG_pv$. Recall that
   $T_p\Sigma$ is a section for the slice representation at $p$ (see \cite{pt}, Theorem 4.6).
   Therefore $T_vG_pv$ is a subspace of $T_p^\perp G(p)$ orthogonal to $T_p\Sigma$ with
\be \label{eq:dim} \dim T_vG_pv=\dim T_p^\perp G(p)-\dim \Sigma.\ee Moreover, again by
the $G$-invariance of $M^n$,  we have that $G(p)\subset M^n$, and
hence $T_pG(p)\subset T_pM^n$. Using (\ref{eq:dim}), we conclude that
$$\begin{array}{l}\dim T_pM^n\geq \dim G(p)+\dim \Sigma +\dim T_vG_pv\vspace{1ex}\\
\hspace*{10.5ex}=\dim G(p) +\dim T_p^\perp G(p)=n+1,\end{array}$$
a contradiction.

  In order to prove the direct statement, it suffices to show that at each point $p\in L$
  which is a singular point of the $G$-action the subset
$$H=\{\gamma'(0)\,|\,\gamma\colon\,(-\epsilon,\epsilon)\to \R^{n+1}\,\,\mbox{is a curve in $\R^{n+1}$ with $\gamma(-\epsilon,\epsilon)\subset M^n$ and }\,\, \gamma(0)=p\}$$
is an $n$-dimensional subspace of $\R^{n+1}$, the tangent space of $M^n$ at $p$.
Clearly, we have
$$H=T_pG(p)\bigcup_{g\in G_p}gT_pL.$$
We use again that $T_p\Sigma$ is a section of the slice representation at $p$ and that,
in addition, the Weyl group for the slice representation
is $W(\Sigma)_p=W\cap G_p$ (\cite{pt}, Theorem 4.6).  Since $L$ is invariant under
$W$ by assumption, it follows that  $T_pL$ is invariant under
$W(\Sigma)_p$. Let $\xi\in \Sigma$ be a unit vector normal to $L$ at $p$. Then,
for any principal vector  $v\in T_pL\subset T_p\Sigma$ of the slice representation
at $p$ it follows from Proposition \ref{prop:weyl} that $G_pv$ lies in the affine
hyperplane ${\cal \pi}$
of $\R^N$ through $p$ orthogonal to $\xi$. Therefore $gT_pL\subset {\cal \pi}$
for every $g\in G_p$.  Since $T_pG(p)$ is orthogonal to $\Sigma$, we conclude that
$H\subset {\cal \pi}$, and hence $H= {\cal \pi}$ by dimension reasons.\qed
\begin{remark}{\em  Theorems  \ref{thm:main} and \ref{thm:weyl} yield as a
special case the main theorem in \cite{mpst}: any compact
hypersurface of $\R^{n+1}$ with cohomogeneity one under the action
of a closed connected subgroup of its isometry group is given as
$G(\gamma)$, where $G\subset SO(n+1)$ acts on $\R^{n+1}$ with
cohomogeneity two (hence polarly) and $\gamma$ is a smooth curve
in a (two-dimensional) section $\Sigma$ which is invariant under
the Weyl group $W$ of the $G$-action. We take the opportunity to
point out that starting with a smooth curve $\beta$ in a Weyl
chamber $\sigma$ of $\Sigma$ (which is identified with the orbit
space of the $G$-action) which is orthogonal to the boundary
$\partial \sigma$ of $\sigma$ is not enough to ensure smoothness
of $\gamma=W(\beta)$, or equivalently, of $G(\beta)$, as claimed
in \cite{mpst}. One should require, in addition, that after
expressing $\beta$ locally as a graph $z=f(x)$ over its tangent
line at a point of intersection with $\partial \sigma$, the
function $f$ be  {\em even\/}, that is, all of its derivatives of
odd order vanish and not only the first one.}
\end{remark}

\section[Proofs of   Corollaries \ref{cor:rev} and  \ref{cor:nneg}.]{Proofs of
Corollaries  \ref{cor:rev} and
\ref{cor:nneg}. }

For the proof of Corollary \ref{cor:rev} we need another known
property of polar representations, a simple proof of which is included for the sake
of completeness

\begin{lemma} \label{le:pol2} Let $\tilde{G}\subset SO(n+1)$ act polarly and have a principal
orbit $\tilde{G}(q)$ that is not
full, i.e.,  the linear span $V$ of $\tilde{G}(q)$ is a proper
subspace of $\R^{n+1}$. Then $\tilde{G}$ acts trivially on
$V^{\perp}$.
\end{lemma}
\proof Let $v\in V^\perp$. Then $v$ belongs to the normal spaces
of $\tilde{G}(q)$ at any point, hence admits a unique extension to
an equivariant normal vector field $\xi$ along $\tilde{G}(q)$.
Moreover, the shape operator $A_\xi$ of $\tilde{G}(q)$ with
respect to $\xi$ is identically zero, for $A_v$ is clearly zero
and $\xi$ is equivariant. Now, since $\tilde{G}$ acts polarly, the
vector field $\xi$ is parallel in the normal connection. It
follows that $\xi$ is a constant vector in $\R^{n+1}$, which means
that $\tilde{G}$ fixes $v$.\vspace{2ex}\qed

 \noindent {\em Proof of Corollary
\ref{cor:rev}.} We know from Theorem \ref{cor:polar} that there
exists an orthogonal representation $\Psi\colon\,G\to SO(n+1)$
such that $\tilde{G}=\Psi(G)$ acts polarly on $\R^{n+1}$ and
$f\circ g=\Psi(g)\circ f$ for every $g\in G$. We claim that any of
the conditions in the statement implies that $f$ immerses some
$G$-principal orbit $Gp$ in $\R^{n+1}$ as a round sphere (circle,
if $k=n-1$).   Assuming the claim, it follows from Lemma \ref{le:pol2} that $\tilde{G}$
fixes the orthogonal complement $V^\perp$ of the linear span $V$ of $\tilde{G}p$, hence
$f$ is a rotation hypersurface with $V^\perp$ as axis. Moreover,
if $k\neq n-1$ then $f|_{Gp}$ must be  an embedding by a standard
covering map argument. From this and $f\circ g=\Psi(g)\circ f$ for
every $g\in G$ it follows that any $g$ in the kernel of $\Psi$
must fix any point of $Gp$. Since $Gp$ is a principal orbit, this
easily implies that $g=\mbox{id}$. Therefore $\Psi$ is an
isomorphism of $G$ onto $\tilde{G}$.

   We now prove the claim. By Proposition  \ref{prop:pt1} we have that
$f$ immerses $Gp$ as an isoparametric submanifold. If condition
$(i)$ holds  then the first normal spaces of $f|_{Gp}$ in
$\R^{n+1}$ are one-dimensional. By Theorem \ref{thm:s}-$(ii)$ and
a well-known result on reduction of codimension of isometric
immersions (cf.  \cite{daj}, Proposition $4.1$), we obtain that
$f(Gp)=\tilde{G}(f(p))$ is contained as a hypersurface in some
affine hyperplane ${\cal H}\subset \R^{n+1}$, and hence $f(Gp)$
must be a round hypersphere of ${\cal H}$ (a circle if $k=n-1$).
Moreover, condition $(i)$ is automatic if $k=n-1$, for in this
case the principal orbit of maximal length must be a geodesic.
Now, if $(iii)$ holds, let $Gp$ be a principal orbit  such that
the position vector of $f$ is nowhere tangent to $f(M^n)$ along
$Gp$. Then any normal vector $\tilde{\xi}$ to $f|_{Gp}$ at $gp\in
Gp$ can be written as $\tilde{\xi}=af(gp)+f_*(gp)\xi_{gp}$, with
$\xi_{gp}$ normal to $Gp$ in $M^n$. Therefore the shape operator
$A^{f|_{Gp}}_{\tilde{\xi}}$ is a multiple of the identity tensor,
hence $f|_{Gp}$ is umbilical. Since, by $(ii)$, the dimension of
$Gp$ can be assumed to be at least two, the claim is proved also
in this case. As for conditions $(iv)$ and $(v)$, the claim
follows from Proposition \ref{prop:isop1}-$(i)$ and the last
assertion in Proposition \ref{prop:isop1}-$(ii)$, respectively.
\vspace{2ex}\qed

\noindent {\em Proof of Corollary \ref{cor:nneg}.} By Proposition
\ref{prop:pt1} and part $(ii)$ of Proposition \ref{prop:isop1} we
have that the orbit $\tilde{G}(f(p))=f(Gp)$ of $\tilde{G}$ is an
extrinsic product $\Sf_1^{n_1-1}\times \ldots\times \Sf_k^{n_k-1}$
of round spheres or circles. In particular, this implies that the
orthogonal decomposition $\R^{n+1}=\bigoplus_{i=0}^k \R^{n_i}$,
where $\R^{n_i}$ is the linear span of $\Sf_i^{n_i-1}$ for $1\leq
i\leq k$, is $\tilde{G}$-stable, that is,  $\R^{n_i}$ is
$\tilde{G}$-invariant for $0\leq i\leq k$ (cf. \cite{got}, Lemma
$6.2$). By  \cite{d}, Theorem $4$  there exist connected Lie
subgroups $G_1,\ldots, G_k$ of $\tilde{G}$ such that
  $G_i$ acts on $\R^{n_i}$ and the action of $\bar{G}=G_1\times\ldots\times G_k$
 on $\R^{n+1}$ given by
 $$ (g_1\ldots
 g_k)(v_0,v_1,\ldots,v_k)=(v_0,g_1v_1,\ldots,g_kv_k)$$
 is orbit equivalent to the action of $\tilde{G}$.
Moreover, writing $q=f(p)=(q_0,\ldots, q_k)$ then
 $\tilde{G}(q)=\{q_0\}\times G_1(q_1)\times\ldots \times
 G_k(q_k)$, and hence $G_i(q_i)$ is a hypersphere of $\R^{n_i}$ for $1\leq i\leq k$.
  The last assertion is clear.\vspace{2ex}\qed

\section[Proof of Theorem \ref{thm:warp}.]{Proof of  Theorem \ref{thm:warp}. }

 Let
$L^k\times_{\rho}N^{n-k}$ be a warped product with  $N^{n-k}$
connected and complete and let
$\psi\colon\,L^k\times_{\rho}N^{n-k}\to U$ be an isometry onto an
open dense subset $U\subset M^n$.  Since $M^n$ is a compact
Riemannian manifold isometrically immersed in Euclidean space as a
hypersurface, there exists an open subset $W\subset M^n$ with
strictly positive sectional curvatures. The subset $U$ being open
and dense in $M^n$,   $W\cap U$ is a nonempty open set. Let
$L_1\times N_1$ be a  connected open subset of $L^k\times N^{n-k}$
that is mapped into $W\cap U$ by $\psi$.
 Then the sectional curvatures of $L^k\times_{\rho}N^{n-k}$ are strictly
 positive on $L_1\times N_1$.

 For a fixed  $x\in L_1$, choose a unit vector $X_x\in
 T_xL_1$. For each $y\in N^{n-k}$, let $X_{(x,y)}$ be the unique unit
 horizontal vector in $T_{(x,y)}(L_1\times N^{n-k})$ that projects onto
 $X_x$ by $(\pi_1)_*(x,y)$. Then the sectional curvature of
 $L^k\times_{\rho}N^{n-k}$ along a plane $\sigma$ spanned by $X_{(x,y)}$ and any
 unit  vertical vector $Z_{(x,y)}\in T_{(x,y)}(L_1\times N^{n-k})$ is
 given by $$K(\sigma)=-\hess \rho(x)(X_x,X_x)/\rho(x).$$ Observe that $K(\sigma)$
 depends neither on $y$
 nor on the vector $Z_{(x,y)}$.
 Since $K(\sigma)>0$ if $y\in N_1$,  the same holds for any
 $y\in N^{n-k}$. In particular, $L_1\times_{\rho} N^{n-k}$ is free of flat
 points.

   If $n-k\geq 2$, it follows from  \cite{dt}, Theorem $16$  that $f\circ \psi$ immerses
   $L_1\times_{\rho} N^{n-k}$ either as a
  rotation hypersurface or as the extrinsic product of an Euclidean factor $\R^{k-1}$ with
  a cone over a hypersurface of $\Sf^{n-k+1}$. The latter possibility is ruled out by the fact
   that the sectional
curvatures of $L^k\times_{\rho}N^{n-k}$ are strictly positive on
$L_1\times N_1$. Thus the first possibility holds, and in
particular $f\circ \psi$ immerses each leaf $\{x\}\times N^{n-k}$,
$x\in L_1$, isometrically onto a round $(n-k)$-dimensional sphere.
It follows that $N^{n-k}$ is isometric to a round sphere.

     In any case, $\mbox{Iso}^0(N^{n-k})$ acts transitively on $N^{n-k}$ and each
     $g\in \mbox{Iso}^0(N^{n-k})$
     induces an isometry $\bar{g}$ of $L^k\times_{\rho}N^{n-k}$
     by defining
     $$\bar{g}(x,y)=(x,g(y)),\,\,\mbox{for all }\,\,(x,y)\in L^k\times
     N^{n-k}.$$
     The map $g\in \mbox{Iso}^0(N^{n-k})\mapsto \bar{g}\in \mbox{Iso}(L^k\times_{\rho}N^{n-k})$
      being clearly continuous, its image $\bar{G}$  is a closed connected
     subgroup of $\mbox{Iso}(L^k\times_{\rho}N^{n-k})$. For each
     $\bar{g}\in \bar{G}$, the induced isometry $\psi\circ \bar{g}\circ \psi^{-1}$ on $U$
     extends uniquely to an isometry of $M^n$. The orbits of
     the induced action of $\bar{G}$ on $U$ are the images by $\psi$ of the leaves
     $\{x\}\times N^{n-k}$, $x\in L^k$,
     hence are umbilical in $M^n$. Moreover, the
normal spaces to the (principal) orbits of $\bar{G}$ on $U$ are
the images by $\psi_*$ of the horizontal subspaces  of
$L^k\times_{\rho}N^{n-k}$. Therefore, they define an integrable
distribution  on $U$, whence on the whole regular part of $M^n$.
Thus, the action of $\bar{G}$ on $M^n$ is locally polar with
umbilical principal orbits. We conclude from Corollary
\ref{cor:rev}-$(iii)$ that $f$ is a rotation hypersurface.\qed

\vspace{.2in} {\renewcommand{\baselinestretch}{1}

\hspace*{-20ex}\begin{tabbing} \indent     \= Universidade Federal
Fluminense
\indent\indent  \=  Universidade Federal de S\~{a}o Carlos \\
\> 24020-140 -- Niteroi -- Brazil  \>
13565-905 -- S\~{a}o Carlos -- Brazil \\
\> E-mail: ion.moutinho@ig.com.br \> E-mail: tojeiro@dm.ufscar.br
\end{tabbing}}

\end{document}